\documentclass[a4paper,12pt]{article}
\usepackage{amssymb,times,amsmath,amsthm} 
\usepackage{hyperref}
\hypersetup{pdfstartview=FitH}
\hypersetup{pdfpagemode=FitWidth} 
\hypersetup{bookmarksnumbered}
\newcommand{\n}{\newcommand}
\n{\dd}{\DeclareMathOperator}\n{\oo}{\operatorname}
\n{\aeq}{\begin{equation}}\n{\zeq}{\end{equation}}
\n{\am}{\begin{pmatrix}}\n{\zm}{\end{pmatrix}}
\n{\mb}{\mathbb}
\n{\bb}{\bigskip}
\newcommand{\bi}{\binom}
\n{\cl}{\centerline} 
\n{\DD}{\Delta}
\n{\de}{euclidean division }
\n{\dl}{jet}
\n{\dls}{jets}
\dd{\DL}{J}
\n{\dr}{\partial}
\n{\ds}{\displaystyle}
\n{\e}{\equiv}
\n{\edo}{differential equation}
\n{\E}{$E$-euclidean}
\n{\el}{element}
\n{\els}{elements}
\n{\f}{\varphi} 
\n{\fr}{rational fraction}
\n{\frs}{rational fractions}
\n{\frd}{rational fraction defined at $a$}
\n{\frds}{rational fractions defined at $a$}
\n{\frg}{generalized rational fraction}
\n{\frgs}{generalized rational fractions}
\n{\ii}{\hskip5mm\relax}
\n{\inv}{{-1}}
\dd{\Ker}{Ker}
\n{\nc}{complex number}
\n{\p}{polynomial}%\n{\ps}{polynômes}
\n{\pg}{generalized polynomial}
\n{\pgs}{generalized polynomials}
\n{\pgcv}{variable coefficient generalized polynomial}
\n{\pgcvs}{variable coefficient generalized polynomials}
\n{\scr}{\scriptstyle}\n{\deux}[2]{\stackrel{\scr #1}{#2}}
\n{\then}{\Rightarrow}
\n{\spd}{sums, products and derivatives}
\n{\sel}{Laurent series in $X-a$}
\n{\sels}{Laurent series in $X-a$}
\n{\ti}{\times}
\n{\w}{\wedge}
\n{\C}{\mb{C}}\n{\R}{\mb{R}}\n{\N}{\mb{N}}\n{\Q}{\mb{Q}}\n{\Z}{\mb{Z}}
\newtheorem{ttt}{Theorem}\n{\att}
{\begin{ttt}}\n{\ztt}{\end{ttt}}
\newtheorem{cor}[ttt]{Corollary} \n{\acc}{\begin{cor}}\n{\zcc}{\end{cor}}
\newtheorem{lem}[ttt]{Lemma} \n{\al}{\begin{lem}}\n{\zl}{\end{lem}}
\newtheorem{prop}[ttt]{Proposition} \n{\ap}{\begin{prop}}\n{\zp}{\end{prop}}
\theoremstyle{definition}%\newtheorem{thm}{Théorème} 
\newtheorem{xx}[ttt]{Exercise} \n{\ax}{\begin{xx}}\n{\zx}{\end{xx}}
\newtheorem{nota}[ttt]{Convention} \n{\an}{\begin{nota}}\n{\zn}{\end{nota}}
\newtheorem{notas}[ttt]{Notation} \n{\ans}{\begin{notas}}\n{\zns}{\end{notas}}
\newtheorem{rmq}[ttt]{Remark} \n{\ar}{\begin{rmq}}\n{\zr}{\end{rmq}}
\begin{document}
{\footnotesize For the last version of this text, type {\em didrygaillard} on Google. Date of this version: Tue Dec 23 09:13:38 CET 2008. Jean-Marie Didry and Pierre-Yves Gaillard} \vfill

\cl{\LARGE Around the Chinese Remainder Theorem}\vfill

\tableofcontents\vfill\newpage

\section{Introduction}
The purpose of this text is to provide a convenient formulation of some folklore results. The five main statements are 
Theorem \ref{rp1} page \pageref{rp1}, 
Theorem \ref{rp2} page \pageref{rp2}, 
Theorem \ref{rp3} page \pageref{rp3}, 
Theorem \ref{rp4} page \pageref{rp4}, and 
Theorem \ref{rp5} page \pageref{rp5}. \bb

\ii By ``\p'' we mean ``complex coefficient \p\ in the indeterminate $X$''. \bb

\ii We show that

\begin{itemize}
\item the computation of the quotient and the remainder of the \de\ of a \p\ by a nonzero \p,
\item the partial fraction decomposition of a \fr, 
\item the computation of an inductive sequence,
\item the exponentiation of a matrix,
\item the integration of an order $n$ constant coefficient linear \edo,
\end{itemize}

follow from a unique, simple and obvious formula: Formula (\ref{tg}) page~\pageref{tg}, which we call {\em Taylor-Gauss Formula}. \bb

\ii The field $\C(X)$ of \frs\ and the ring $E$ of entire functions are two important examples of differential rings containing $\C[X]$. The subring $\C[X]$ ``controls'' $E$ in the sense that an entire function can be euclideanly divided by a nonzero \p, the remainder being a \p\ of degree strictly less than that of the divisor. Among the rings having this property, one contains all the others: it is the differential ring 
\begin{equation}\label{ce}
\prod_{a\in\C}\C[[X-a]].
\end{equation}
Any torsion $\C[X]$-module is a module over the ring (\ref{ce}), and this ring is universal for this property. In particular any element $f$ of (\ref{ce}) can be evaluated on a square matrix $A$, the matrix $f(A)$ being by definition $R(A)$, where $R$ is the remainder of the \de\ of $f$ by a nonzero \p\ annihilating $A$. Since there is an obvious formula for this remainder (the Taylor-Gauss Formula), we are done. An important example consists in taking for $f$ the exponential function, viewed as  the element of (\ref{ce}) whose $a$-th component is the Taylor series of $e^X$ at $a$. We recover of course the usual notion of exponential of a matrix, but cleared form its artificial complications. \bb

\ii It is handy to embed $\C(X)$ and $E$ into the differential ring $$\prod_{a\in\C}\C((X-a)),$$ which we use as a huge container. A side advantage of this ring is that it shortcuts the usual (particularly unilluminating) construction of the field of \frs\ (as a differential field) from the ring of \p s. \bb

%\ii The main idea used in this text shows up in {\em Local Class Field Theory} by Serre (see [1]). 
%
\section{Laurent Series}
Let $a$ be a \nc. A {\bf\sel} is an expression of the form 
$$f=f(X)=\sum_{n\in\Z}f_{a,n}\ (X-a)^n,$$
where $(f_{a,n})_{n\in\Z}$ is a family of complex numbers for which there is an integer $n_a$ such that $n<n_a$ implies $f_{a,n}=0$. \bb

\ii We define the operations of addition, multiplication and differentiation on the \sels\ by
$$(f+g)_{a,n}=f_{a,n}+g_{a,n},$$
$$(f\,g)_{a,n}=\sum_{p+q=n}f_{a,p}\ g_{a,q},$$
$$(f')_{a,n}=(n+1)\ f_{a,n+1},$$
and we check that these operations have the same properties as on \p s.
\att Let $f$ be a \sel. If $f\not=0$, then there is a unique Laurent series $g$ in $X-a$ such that $f\,g=1$.\ztt
{\bf Proof.} Exercise.\bb

\ii Put $g=1/f=\frac{1}{f}$ and $h\,g=h/f=\frac{h}{f}$ if $h$ is a \sel.
\section{Generalized Rational Fractions}
A {\bf \frg} is a family $f=(f_a)_{a\in\C}$ each of whose members $f_a$ is a \sel. The complex numbers $f_{a,n}$ are called the {\bf coefficients} of $f$. The \frgs\ are added, multiplied and differentiated componentwise. \bb

\ii To each \p\ $P$ is attached the \frg
\begin{equation}\label{P}
\left(\sum_{n=0}^\infty\frac{P^{(n)}(a)}{n!}\ 
(X-a)^n\right)_{a\in\C}.
\end{equation}
As the formal sum in the parenthesis contains only a finite number of nonzero terms, it can be viewed as a \p. As such, it is of course equal to the \p\ $P$. Hence the \spd\ of \p s as \p s coincide with their \spd\ as \frgs. These facts prompt us to designate again by $P$ the \frg\ (\ref{P}).\bb

\ii We can now define a {\bf \fr} as being a \frg\ obtained by dividing a \p\ by a nonzero \p. The \spd\ of \frs\ are \frs. If a nonzero \frg\ $f=(f_a)_{a\in\C}$ is a \fr, then $f_a$ is nonzero for all $a$. \bb

\ii For all \frg\ $f$ and all \nc\ $a$ put
$$\mu(a,f):=\inf\ \{n\in\Z\ |\ f_{a,n}\not=0\}$$
with the convention $\inf\ \varnothing=+\infty$, and say that $\mu(a,f)$ is the {\bf multiplicity} of $a$ as zero, or root, of $f$. We have 
$$\mu(a,f g)=\mu(a,f)+\mu(a,g).$$ 
If $\mu(a,f)\ge0$ say that $f$ is {\bf defined at} $a$, and denote $f_{a,0}$ by $f(a)$. If a \frg\ $f$ is defined at a \nc\ $a$, then we have
$$f_a=\sum_{n=0}^\infty\frac{f^{(n)}(a)}{n!}\ (X-a)^n.$$
Any sum, product or derivative of \frgs\ defined at $a$ is a \frg\ defined at $a$. \bb

\ii For all \frg\ $f$, all \nc\ $a$, and all integer $\mu$ set 
$$\DL_a^\mu(f):=\sum_{n\le\mu}f_{a,n}\ (X-a)^n,$$
and say that this \sel\ is the {\bf order} $\mu$ {\bf \dl\ of $f$ at $a$}. If $f$ and $g$ are \frgs\ defined at $a$, we have
$$\DL_a^\mu(f+g)=\DL_a^\mu(f)+\DL_a^\mu(g),\quad
\DL_a^\mu(f\,g)
=\DL_a^\mu\Big(\DL_a^\mu(f)\DL_a^\mu(g)\Big).$$

\ii Let $a$ be a \nc, let $\mu$ be an integer, and let $f$ and $g$ be two \frgs. 
\ax\label{mod} Show that the following conditions are equivalent
\begin{enumerate}
\item $\mu(a,f-g)\ge\mu$,
\item $(X-a)^{-\mu}\ (f-g)$ is defined at $a$,
\item $\DL_a^{\mu-1}(f)=\DL_a^{\mu-1}(g)$.
\end{enumerate}
\zx 

\ii If these conditions are satisfied and if $f$ and $g$ are defined at $a$, then we say that $f$ and $g$ are {\bf congruent modulo} $(X-a)^\mu$ and we write $$f\e g\bmod(X-a)^\mu.$$
We thus have
$$\DL_a^{\mu-1}(f)\e f\bmod(X-a)^\mu,$$
as well as
$$\left.
\begin{array}{c}
 f_1\e g_1\bmod (X-a)^\mu\\ [1em]
 f_2\e g_2\bmod (X-a)^\mu
\end{array}\right\}\then\left\{
\begin{array}{c}
 f_1+f_2\e g_1+g_2\bmod (X-a)^\mu\\  [1em]
 f_1\,f_2\e g_1\,g_2\bmod (X-a)^\mu.
\end{array}\right.$$
\ax\label{pr} Assume $f$ is defined at $a$ and $\mu\ge0$. Let $R$ be a degree $<\mu$ \p. Show that the following conditions are equivalent. 
\begin{enumerate}
\item $R=\DL_a^{\mu-1}(f)$,
\item $R\e f\bmod(X-a)^\mu$,
\item $(X-a)^{-\mu}\ (f-R)$ is defined at $a$. 
\end{enumerate}
In other words 
$$\left(\frac{f-\DL_a^{\mu-1}(f)}{(X-a)^\mu}\ ,\ 
\DL_a^{\mu-1}(f)\right)$$ 
is the unique pair $(q,R)$ such that
\begin{enumerate}
\item $q$ is a \frg\ defined at $a$, 
\item $R$ is a \p\ of degree $<\mu$, 
\item $f(X)=(X-a)^\mu\,q(X)+R(X)$.
\end{enumerate}\zx

\ii If we wish to extend the result of Exercise \ref{pr} to the division by an arbitrary \p\ $D$, it is natural to restrict to the \frgs\ defined at all point of $\C$. 
We thus define a {\bf \pg} as being a \frg\ defined at all point of $\C$. Any sum, product or derivative of \pgs\ is a \pg. 
\ans Fix a nonconstant \p\ $D$ \label{D} and, for all complex number $a$, let $\mu_a$ be the multiplicity of $a$ as a root of $D$. 
\zns

\ax\label{pg}Let $f$ be a \pg. Show that the \frg\ $f/D$ is a \pg\ if and only if $f\e0\bmod(X-a)^{\mu_a}$ for all $a$. 
\zx

\ii We admit 
\att[Fundamental Theorem of Algebra] There a nonzero  complex number $c$ satisfying
$$D(X)=c\prod_{D(a)=0}(X-a)^{\mu_a}.$$
\ztt
\ax\label{p}Let $P$ be a \p. Show that the following conditions are equivalent
\begin{enumerate}
\item $P\e0\bmod(X-a)^{\mu_a}$ for all $a$, 
\item $P/D$ is a \pg,
\item $P/D$ is a \p,
\item $P/D$ is a degree $\deg P-\deg D$ \p.
\end{enumerate} 
\zx
\section{Chinese Remainder Theorem}
\att[Chinese Remainder Theorem] For all \pg\ $f$ and all nonconstant \p\ $D$ there is a unique couple $(q,R)$ such that\label{tc}
\begin{enumerate}
\item $q$ is a \pg, 
\item $R$ is a \p\ of degree $<\deg D$, 
\item $f=D\,q+R$.
\end{enumerate}
In Notation \ref{D} page \pageref{D} we have the {\bf Taylor-Gauss Formula} (see [1])
\aeq\label{tg}\boxed{\boxed{
R(X)=\sum_{D(a)=0}\ \DL_a^{\mu_a-1}\!\!\left(f(X)\
\frac{(X-a)^{\mu_a}}{D(X)}\right)
\frac{D(X)}{(X-a)^{\mu_a}}}}\zeq\ztt

\ii Say that $R$ is the {\bf remainder of the \de} of $f$ by $D$, and that $q$ is its {\bf quotient}. \bb

{\bf Proof.} Uniqueness. Assume
$$f=D\ q_1+R_1=D\ q_2+R_2$$
(obvious notation). We get
$$\frac{R_2-R_1}{D}=q_1-q_2$$
and Exercise \ref{p} (2. $\then$ 4.) implies $R_1=R_2$ and thus $q_1=q_2$. \bb

\ii Existence. Put
$$s_{(a)}(X):=\DL_a^{\mu_a-1}\!\!\left(f(X)\ 
\frac{(X-a)^{\mu_a}}{D(X)}\right)
\frac{1}{(X-a)^{\mu_a}}\quad,\quad s:=\sum\ s_{(a)},$$
and note that $\frac{f}{D}-s_{(a)}$ is defined at $a$ by Exercise~\ref{pr}, that $\frac{f}{D}-s$ is a \pg\ $q$, and that $s D$ is a degree $<\deg D$ \p\ $R$. QED\bb

{\bf Example.} Let $a$ and $b$ be two complex numbers. For $n\ge0$ write $s_n$ for the sum of the degree $n$ monomials in $a$ and $b$. The remainder of the \de\ of the \p\ $\sum_{n\ge0}a_n\ X^n$ by $(X-a)(X-b)$ is
$$\sum_{n\ge1}\ a_n\ s_{n-1}\ X+a_0
-a\ b\ \sum_{n\ge2}\ a_n\ s_{n-2}.$$
\acc {\em (First Serret Formula --- see [1].)} The principal part of $g:=f/D$ at $a$ is \label{s1}
$$\boxed{\boxed{
\DL_a^{\mu_a-1}\!\!\Big(g(X)(X-a)^{\mu_a}\Big)
(X-a)^{-\mu_a}}}
$$\zcc

\ii There are partial analogues to Theorem \ref{tc} and Corollary \ref{s1} over an arbitrary commutative ring (see Section~\ref{aq}), but the uniqueness of the partial fraction decomposition disappears. For instance we have, in a product of two nonzero rings, 
$$\frac{1}{X-1}-\frac{1}{X}
=\frac{(1,-1)}{X-(1,0)}+\frac{(-1,1)}{X-(0,1)}\quad.$$
\att Let $P$ be a \p, let $f$ be the \fr\ $P/D$, let $Q$ be the quotient, which we assume to be nonzero, of the \de\ of $P$ by $D$, and let $q$ be the degree of $Q$. Then $Q$ is given by the {\bf Second Serret Formula} (see [1])~: 
$$\boxed{Q(X^{-1})
=\DL_0^q(f(X^{-1})\ X^q)\ X^{-q}}$$
\ztt
{\bf Proof.} If $d$ is the degree of $D$ and if $R$ is the remainder of the \de\ of $P$ by $D$, then the \fr
$$g(X):=X^\inv\Big(f(X^\inv)-Q(X^\inv)\Big)
=\frac{X^{d-1}\ R(X^\inv)}{X^d\ D(X^\inv)}$$
is defined at 0, and we have
$$X^q\ f(X^\inv)-X^q\ Q(X^\inv)=X^{q+1}\ g(X)
\e0\bmod X^{q+1},$$
which, as $X^q\,Q(X^\inv)$ is a \p\ of degree at most $q$, implies
$$\DL_0^q\Big(X^q\ f(X^\inv)\Big)-X^q\ Q(X^\inv)=0.$$
QED\bb

\ii In Notation \ref{D} page \pageref{D} fix a root $a$ of $D$, denote by $B$ the set of all other roots, and, for any map $u:B\to\N,b\mapsto u_b$, denote by $|u|$ the sum of the $u_b$. 
\att The coefficient $c_{a,k}$ of $(X-a)^k$ in 
$\DL_a^{\mu_a-1}(\frac{(X-a)^{\mu_a}}{D(X)})$ is \label{fe}
\begin{equation}\label{efe}
c_{a,k}=(-1)^k\sum_{\deux{u\in\N^{B}}{|u|=k}}\ \prod_{b\in B}\ 
\binom{\mu_b-1+u_b}{\mu_b-1}\ \frac{1}{(a-b)^{\mu_b+u_b}}\quad.
\end{equation}
\ztt
{\bf Proof.} It suffices to multiply the \dls
$$\DL_a^{\mu_a-1}\!\!\left(\frac{1}{(X-b)^{\mu_b}}\right)
=\sum_{n=0}^{\mu_a-1}\ \binom{\mu_b-1+n}{\mu_b-1}\ 
\frac{(-1)^n\ (X-a)^n}{(a-b)^{\mu_b+n}}\quad.$$
QED
\acc The \p\ $R(X)$ in the Chinese Remainder Theorem (Theorem~\ref{tc} page \pageref{tc}) is \label{e}
$$R(X)=\ds\sum_{\deux{D(a)=0}{k+n<\mu_a}}\ c_{a,k}\ 
\frac{f^{(n)}(a)}{n!}\ (X-a)^{k+n},$$
where $c_{a,k}$ is given by (\ref{efe}). 
\zcc
\acc The coefficients of the partial fraction decomposition of a \fr\ are integral coefficient \p s in \label{cu} \bb
 
\ii $\bullet$ the coefficients of the numerator, 

\ii $\bullet$ the roots of the denominator, 
 
\ii $\bullet$ the inverses of the differences of the roots of the denominator. \bb

These \p s are homogeneous of degree one in the coefficients of the numerator and depend only on the multiplicities of the roots of the denominator. (We assume the denominator is monic.) 
\zcc 
\section{Exponential}
The most important example of \pg\ is perhaps the exponential
$$e^X:=\left(e^a\sum_{n=0}^\infty\frac{(X-a)^n}{n!}
\right)_{a\in\C},$$
which satisfies
$$\frac{d}{d X}\ e^X=e^X$$
as a \pg. \bb

\ii More generally we can, for any real number $t$, define the \pg
\begin{equation}\label{etx}
e^{t X}
:=\left(e^{a t}\sum_{n=0}^\infty\frac{t^n\ (X-a)^n}{n!}
\right)_{a\in\C},
\end{equation}
and observe the identity between \pgs
$$e^{t X}\ e^{u X}=e^{(t+u) X}.$$

\ii We would like to differentiate $e^{t X}$ not only with respect to $X$ but also with respect to $t$. To do this properly we must allow the coefficients of a \pg\ to be not only complex constants, but $C^\infty$ functions from $\R$ to $\C$. Let us call {\bf \pgcv} a \pg\ whose coefficients are $C^\infty$ functions from $\R$ to $\C$. We then have the following identity between \pgcvs
$$\frac{\partial e^{t X}}{\partial t}=X\,e^{t X}.$$

\ii One is often lead to compute \dls\ of the form $$\DL_a^\mu(e^{t X}f(X)),$$ where $f(X)$ is a \frd. This \dl\ is the unique degree $\le \mu$ \p\ satisfying
\begin{equation}\label{w1}
\DL_a^\mu(e^{t X}f(X))\e e^{a t}\ \DL_a^\mu (f(X))\
\sum_{n=0}^\mu\frac{t^n\ (X-a)^n}{n!}\ \bmod(X-a)^{\mu+1}.
\end{equation}

\ii Let us sum this up by 
\att In Notation \ref{D} page \pageref{D} the remainder of the \de\ of $e^{t X}$ by $D(X)$ is given by the {\bf Wedderburn Formula} (see [1])
\begin{equation}\label{w2}
\boxed{\sum_{D(a)=0}\ \DL_a^{\mu_a-1}\!\!\left(e^{t X}\
\frac{(X-a)^{\mu_a}}{D(X)}\right)
\frac{D(X)}{(X-a)^{\mu_a}}}
\end{equation}
the \dl\ being given by (\ref{w1}). 
\ztt
\section{Matrices}
Let $A$ a complex square matrix. 
\ax\label{pm} Show that there is a unique monic \p\ $D$ which annihilates $A$ and which divides any \p\ annihilating $A$. 
\zx 

\ii Say that $D$ is the {\bf minimal \p} of $A$. 
\ax Let $f$ be a \pg, $D$ a \p\ annihilating $A$, and $R$ the remainder of the \de\ of $f$ by $D$. Show that the matrix $R(A)$ does not depend on the choice of the annihilating \p\ $D$. [Suggestion: use Exercises \ref{p} and \ref{pm}.]
\zx 

\ii It is then natural to set $$f(A):=R(A).$$
\ax Let $f,g$ be two \pgs\ and $D$ the minimal \p\ of $A$. Show $f(A)=g(A)\iff D$ divides $f-g$.
\zx 
\ax Show $(f\,g)(A)=f(A)\,g(A)$. \zx 
\att The matrix $e^{t A}:=f(A)$, where $f(X)=e^{t X}$ (see (\ref{etx})), depends differentiably on $t$ and satisfies $\frac{d}{d t}\ e^{t A}=A\ e^{t A},e^{0A}=1$.  
\ztt
\ax Prove the above statement.\zx 
\ax Let $$D(X)=X^q+a_{q-1}\,X^{q-1}+\dots+a_0$$ be a monic \p, $(e_j)$ the canonical basis of $\C^q$, and $A$ the $q$ by $q$ matrix characterised \label{co} by
$$j<q\then A\,e_j=e_{j+1},\quad 
A\,e_q=-a_0\,e_1-\dots-a_{q-1}\,e_q.$$
Show that $D$ is the minimal \p\ of $A$.\zx 

\ii Let $D$ be a degree $q\ge1$ polynomial, let $f$ be a \pg, and $b_{q-1}\,X^{q-1}+\dots+b_0$ the remainder of the \de\ of $f$ by $D$. 
\ax Show $f(A)\,e_1=b_0\,e_1+\dots+b_{q-1}\,e_q$. \label{fA}
\zx  
\section{Inductive Sequences}
Let $D$ be a degree $q\ge1$ monic polynomial; let $\C^\N$ be the set of complex number sequences; let $\DD$ be the shift operator mapping the sequence $u$ in $\C^\N$ to the sequence $\DD u$ in $\C^\N$ defined by $(\DD u)_t=u_{t+1}$; let $f$ be in $\C^\N$; let $c_0,\dots,c_{q-1}$ be in $\C$; let $y$ be the unique element of $\C^\N$ satisfying 
$$D(\DD)\ y=f,\quad y_n=c_n\mbox{ for all } n<q\ ;$$
for $(n,t)$ in $\N^2$ denote by $g_n(t)$ the coefficient of $X^n$ in the remainder of the \de\ of $X^t$ by $D$. 
\att If $t\ge q$ then $\ds y_t=\sum_{n<q}\ c_n\ g_n(t)
+\sum_{k<t}\ g_{q-1}(t-1-k)\ f_k.$
\ztt
{\bf Proof.} Introduce the sequence of vectors $x_t:=(y_t,\dots,y_{t+q-1})$. We have 
$$x_{t+1}=B\,x_t+f_t\,e_q,\quad x_0=c,$$
where $e_q$ is the last vector of the canonical basis of $\C^q$, and $B$ is the transpose of the matrix $A$ in Exercise \ref{co}. Hence
$$x_t=B^t c+f_0\,B^{t-1}\,e_q+f_1\,B^{t-2}\,e_q+\dots +f_{t-1}\,e_q.$$
It suffices then to take the first component of the left and right hand sides, and to invoke Exercise \ref{fA}. QED
\section{Differential Equations}
Let $D$ be a degree $q\ge1$ monic polynomial and $y$ the unique solution to the \edo

\aeq\label{edo}D\left(\frac{d}{d t}\right)y=f(t),\quad
 y^{(n)}(0)=y_n\ \forall\ n<q:=\deg D,\zeq

where $f:\R\to\C$ is a continuous function.
\att We have the {\bf Collet Formula} (see [1])
$$\boxed{y(t)=\sum_{n<q}y_n\ g_n(t)
+\int_0^t g_{q-1}(t-x)\ f(x)\ d x}$$
where $g_n(t)$ is the coefficient of $X^n$ in the remainder of the \de\ of $e^{t X}$ by $D$. 
\ztt

\ii [For $D(X)=X^q$ we recover the Taylor Formula with integral remainder.] \bb 

{\bf Proof.} Putting $v_n:=y^{(n-1)}$, $v_{0n} :=y_{n-1}$ for $1\le n\le q$, and denoting by $e_q$ the last vector of the canonical basis of $\C^q$, Equation (\ref{edo}) takes the form
$$v'(t)-B\ v(t)=f(t)\ e_q,\quad v(0)=v_0,$$
where $B$ is the transpose of the matrix $A$ in Exercise \ref{co}. Applying $e^{-t B}$ we get
$$\frac{d}{d t}\ e^{-t B}\ v(t)=e^{-t B}\ f(t)\ e_q,\quad v(0)=v_0,$$
whence
$$v(t)=e^{t B}\ v_0+\int_0^t f(x)\ e^{(t-x)B}\ e_q\ d x.$$
In view of the Wedderburn Formula (\ref{w2}), it suffices then to take the first component of the left and right hand sides, and to invoke Exercise~\ref{fA}. QED\bb

\ii Let $h$ be a \pgcv, let $f$ and $y$ be two continuous functions from $\R$ to $\C^q$, let $y_0$ be a vector in $\C^q$, and $A$ a $q$ by $q$ complex matrix. If $y$ is differentiable and satisfies
\begin{equation}\label{a}
y'(t)+h(t,A)\ y(t)=f(t),\quad y(0)=y_0,
\end{equation}
then
$$H(t,A):=\int_0^t h(u,A)\ d u\quad\then\quad
\frac{d}{d t}\ e^{H(t,A)}\ y(t)=e^{H(t,A)}\ f(t),$$
whence 
\att The unique solution to (\ref{a}) is given by the {\bf Euler Formula} (see [1])

$$\boxed{y(t)=\exp\left(\int_t^0h(u,A)\ d u\right)y_0
+\int_0^t\exp\left(\int_t^v h(u,A)\ d u\right) f(v)\ d v}$$
\ztt
\section{Euclid}
By ``ring'' we mean in this text ``commutative ring with 1''. For the reader's convenience we prove the Chinese Remainder Theorem. 
\att [Chinese Remainder Theorem] Let $A$ be a ring and $I_1,\dots,I_n$ ideals such that $I_p+I_q=A$ for $p\not=q$. Then the natural morphism \label{tcc} from $A$ to the product of the $A/I_p$ is surjective. Moreover the intersection of the $I_p$ coincides with their product. 
\ztt
{\bf Proof.} Multiplying the equalities $A=I_1+I_p$ for $p=2,\dots,n$ we get
\aeq\label{12n}A=I_1+I_2\cdots I_n.\zeq
In particular there is an $a_1$ in $A$ such that
$$a_1\equiv1\bmod I_1,\quad a_1\equiv0\bmod I_p\ \forall\ p>1.$$
Similarly we can find elements $a_p$ in $A$ such that
$a_p\equiv\delta_{p q}\bmod I_q$ (Kronecker delta). This proves the first claim. Let $I$ be the intersection of the $I_p$. Multiplying (\ref{12n}) by $I$ we get 
$$I=I_1I+II_2\cdots I_n\subset I_1\ (I_2\cap\cdots \cap I_n)\subset I.$$
This gives the second claim, directly for $n=2$, by induction for $n>2$. QED \bb

\ii Let $A$ be a PID (principal ideal domain). Construct the ring $A^\w$ as follows. A family $(a_d)_{d\not=0}$ of \els\ of $A$ indexed by the nonzero \els\ $d$ of $A$ represents an \el\ of $A^\w$ if it satisfies
$$d\ |\ e\then a_d\e a_e\bmod d$$
(where $d\ |\ e$ means ``$d$ divides $e$'') for all pair $(d,e)$ nonzero \els\ of $A$. Two such families $(a_d)_{d\not=0}$ and $(b_d)_{d\not=0}$ represent the same \el\ of $A^\w$ if and only if $$a_d\e b_d\bmod d\quad\forall\ d\not=0.$$ The ring structure is defined in the obvious way. Embed $A$ into $A^\w$ by mapping $a$ in $A$ to the constant family equal to $a$. By abuse of notation we often designate by the same symbol an element of $A^\w$ and one of its representatives. Let $P$ be a representative system of the association classes of prime \els\ of $A$.
\al Let $a=(a_b)_{b\not=0}$ be in $A^\w$ and $d$ a nonzero \el\ of $A$ such that $a_d\e 0\bmod d$. \label{27} Then there is a $q$ in $A^\w$ such that  $a=d\,q$.
\zl
{\bf Proof.} Let $p$ be in $P$ and $i$ the largest integer $j$ such that $p^j$ divides $d$. In other words there is an \el\ $d'$ of $A$ which is prime to $p$ and satisfies $d=p^i d'$. For all nonnegative integer $j$ we have
$$a_{p^{i+j}}\e0\bmod p^i.$$ 
As a result there is an \el\ $a'_j$ of $A$ such that $a_{p^{i+j}}=p^i\,a'_j$. We have
$$p^i\,a'_{j+1}\e a_{p^{i+j+1}}\e a_{p^{i+j}}\e 
p^i\,a'_j\bmod p^{i+j}$$
and thus 
$$a'_{j+1}\e a'_j\bmod p^j.$$
For all $j$ choose $q_{p^j}$ such that 
$$d'\,q_{p^j}\e a'_j\bmod p^j$$
and thus
$$d\,q_{p^j}\e a_{p^j}\bmod p^j.$$
We get
$$d'\,q_{p^{j+1}}\e a'_{j+1}\e a'_j\e d'\,q_{p^j}
\bmod p^j$$
and thus
$$q_{p^{j+1}}\e q_{p^j}\bmod p^j.$$
Let $b$ be in $A$, $b\not=0$, and $P_b$ the (finite) set of those \els\ of $P$ which divide $b$. For $p$ in $P_b$ denote by $i(p)$ the largest integer $j$ such that $p^j$ divides $b$ and choose a solution $q_b$ in $A$ to the congruence system
$$q_b\e q_{p^{i(p)}}\bmod p^{i(p)},\quad p\in P_b,$$
solution which exists by the Chinese Remainder Theorem (Theorem \ref{tcc} page \pageref{tcc}). One then checks that the family $q:=(q_b)_{b\not=0}$ is in $A^\w$ and that we do have $d\,q=a$. QED\bb

\ii See Theorem \ref{rp4} page \pageref{rp4} for a generalization. \bb

\ii Let $a$ be in $A^\w$, let $d$ be a nonzero \el\ of $A$, let $\mu(p)$ be the multiplicity of $p$ in $P$ as a factor of $d$, let 
$$\DL_p^{\mu(p)-1}\!\!\left(a\
\frac{p^{\mu(p)}}{d}\right)$$
be an \el\ of $A$ satisfying
$$\DL_p^{\mu(p)-1}\!\!
\left(a\ \frac{p^{\mu(p)}}{d}\right)
\frac{d}{p^{\mu(p)}}\e a_{p^{\mu(p)}}\bmod p^{\mu(p)},$$
and put 
$$\rho:=\sum_{p\ |\ d}\ \DL_p^{\mu(p)-1}\!\!\left(a\
\frac{p^{\mu(p)}}{d}\right)\frac{d}{p^{\mu(p)}}\in A.$$

\ii Let $B$ be a ring and $A$ an integral domain contained in $B$. Say that $A$ is {\bf principal in} $B$ if for all $b$ in $B$ and all nonzero $d$ in $A$ there is a $q$ in $B$ and an $r$ in $A$ such that $b=d q+r$. \bb
\att[First Main Statement] Let $A$ be a PID, let $a$ be an element of $A^\w$, and let $d$ be a nonzero element of $A$. In the above notation, \label{rp1} the elements $a$, $a_d$ and $\rho$ of $A^\w$ are congruent modulo $d$. Moreover, if $A$ is principal in $B$, then there is a unique $A$-algebra morphism from $B$ to $A^\w$. 
\ztt
{\bf Proof.} Modulo $d$ we have $a\e a_d$ by the Lemma, and $a_d\e\rho$ by the Chinese Remainder Theorem (Theorem \ref{tcc} page \pageref{tcc}). This proves the first claim. Let us check the second one, starting with the uniqueness. Let $f$ be an $A$-linear map from $B$ to $A^\w$. If $b$ and $q$ are in $B$, and if $d\not=0$ and $r$ are in $A$, then the equality $b=d\,q+r$ implies $f(b)=d\,f(q)+r$, and thus $f(b)_d\e r\bmod d$. Consequently there is at most one such map. The existence is proved by putting $f(b)_d:=r$ in the above  notation, and by checking that this formula does define an $A$-linear map from $B$ to $A^\w$. QED \bb

\ii An $A$-module is {\bf torsion} if each of its vectors is annihilated by some nonzero scalar. 
\att If $A$ is principal in $B$, then any torsion $A$-module admits a unique $B$-module structure which extends its \label{w} $A$-module structure. Moreover any $A$-linear map between torsion $A$-modules is $B$-linear. 
\ztt

{\bf Proof.} Let us check the uniqueness. Let $b$ be in $B$~; let $v$ be in our torsion module $V$~; let $d\not=0$ in $A$ satisfy $d v=0$~; let $q$ in $B$ and $r$ in $A$ verify $b=d q+r$. We then get $b v=r v$, which proves the uniqueness. The existence is obtained by setting $b v:=r v$ in the above notation, and by checking that this formula does define a $B$-module structure on $V$ which extends the $A$-module structure. The last assertion is clear. QED \bb

\ii [Let $\oo{Hom}_{c t}(G,H)$ be the group of continuous morphisms from the topological group $G$ to the abelian topological group $H$, and $\widehat{\Z}$ the profinite completion of $\Z$. Equip $\Q/\Z$ with the discrete topology or with the topology induced by that of $\R/\Z$. Equip also $\Q/\Z$ with the $\widehat{\Z}$-module structure provided by Theorem~\ref{w}. We then have  
$$\oo{Hom}_{c t}(\widehat{\Z},\R/\Z)
=\oo{Hom}_{c t}(\widehat{\Z},\Q/\Z)
=\oo{Hom}_{\widehat{\Z}}(\widehat{\Z},\Q/\Z)
=\Q/\Z.$$
We recover the well known fact that $\Q/\Z$ is the dual of $\widehat{\Z}$ in the category of locally compact abelian groups. See Weil's {\em L'int\'egration dans les groupes topologiques}, pp. 108-109.] \bb

\ii Let $L$ be a sublattice of the lattice of ideals of an arbitrary ring $A$. Here are some examples: 
\begin{enumerate}
\item $L$ is the set of nonzero ideals of an integral domain, 
\item $L$ is the set of powers of a fixed ideal of an arbitrary ring, 
\item let $Y$ be a set of prime ideals of $A$, and $L$ the set of those \label{supp} ideals $I$ of $A$ such that any prime ideal containing $I$ is in $Y$, 
\item $L$ is the set of open ideals of a topological ring.
\end{enumerate}
Construct the ring $A^\w$ as follows. An \el\ of $A^\w$ is represented by a family $(a_I)_{I\in L}$ of \els\ of $A$ satisfying
$$I\subset J\then a_I\e a_J\bmod J$$
for all $I,J$ in $L$. Two such families $(a_I)_{I\in L}$ and $(b_I)_{I\in L}$ represent the same \el\ of $A^\w$ if and only if 
$$a_I\e b_I\bmod I\quad\forall\ I\in L.$$ 
The ring structure is defined in the obvious way. By abuse of notation we often designate by the same symbol an element of $A^\w$ and one of its representatives. Map $A$ to $A^\w$ by sending $a$ in $A$ to the constant family equal to $a$. \bb

\ii Call {\bf torsion module} any $A$-module each of whose vector is annihilated by some \el\ of $L$. [In the case of Example \ref{supp} above, an $A$-module is torsion if and only if its support is contained in $Y$ --- Bourbaki, {\em Alg. Com.} II.4.4.] Assume to simplify that the intersection of the \els\ of $L$ reduces to zero, and consider $A$ as a subring of $A^\w$. 
\att[Second Main Statement] Any torsion $A$-module $V$ admits a unique $A^\w$-module structure such that \label{rp2} 
$$a\in A^\w,v\in V,I\in L,Iv=0\then a v=a_I\,v.$$
This $A^\w$-module structure on $V$ extends the $A$-module structure.
\ztt
{\bf Proof.} The uniqueness being obvious, lets us check the existence. For $v$ in $V$ and $I,J$ in $L$ such that $I v=0=J v$, we have $a_I\,v= a_{I+J}\ v$; as this vector depends only on $a$ and $v$, it can be denoted $a v$. Let us show $a\,(v+w)=a v+a w$. If $a$ is in $A^\w$, and if $I,J$ are in $L$ and verify $I v=0=J w$, we have
$$a\ (v+w)=a_{I\cap J}\ (v+w)=a_{I\cap J}\ v+a_{I\cap J}\ w=a\,v+a\,w.$$
QED 
\att[Third Main Statement] Let $B$ a ring containing $A$. Assume that each torsion $A$-module is equipped with a $B$-module structure which extends its $A$-module structure, and that each $A$-linear map between torsion $A$-modules is $B$-linear. \label{rp3} Then there is a unique $A$-algebra morphism $f$ from $B$ to $A^\w$ which satisfies $b\,v=f(b)\,v$ for all vector $v$ of any torsion $A$-module, and all $b$ in $B$. 
\ztt
{\bf Proof.} For all $I$ in $L$ denote by $1_I$ the element 1 of $A/I$. Let us verify the uniqueness. Let $b$ be in $B$ and $I$ in $L$. There is an $a$ in $A$ satisfying
$$a\cdot1_I=b\cdot1_I=f(b)\cdot1_I=f(b)_I\cdot1_I,$$
and thus $f(b)_I\e a\bmod I$. To check the existence, we put $f(b)_I:=a$ in the above notation, and verify that this formula does define an $A$-linear map $f$ from $B$ to $A^\w$ which satisfies $b\,v=f(b)\,v$ for all vector $v$ of any torsion $A$-module, and all $b$ in $B$. [To check that $f(b)_I$ and $f(b)_J$ are congruent modulo $J$ for $I\subset J$ in $L$, we use the assumption that the canonical projection from $A/I$ to $A/J$ is $B$-linear.] QED \bb 

\ii Let $A$ be a Dedekind domain; let $L$ be the set of nonzero ideals of $A$; and let $M$ be the set of maximal ideals of $A$. For any $A$-module $V$, let $V^\w$ be the projective limit of the $V/I V$ with $I$ in $L$; and, for any $P$ in $M$, let $V^\w_P$ be the projective limit of the $V/P^n V$. Then $A^\w$ and $A^\w_P$ are $A$-algebras; by the Chinese Remainder Theorem (Theorem \ref{tcc} page \pageref{tcc}) $A^\w$ is the direct product of the $A^\w_P$; the $A$-module $V^\w$ is an $A^\w$-module; the $A$-module $V^\w_P$ is an $A^\w_P$-module; $V^\w$ is the direct product of the $V^\w_P$. Here is a generalization of Lemma \ref{27} page \pageref{27}. 
\att[Fourth Main Statement] The ring $A^\w$ is $A$-flat. If $V$ is a finitely generated $A$-module, then the natural morphism from $A^\w\otimes_A V$ to $V^\w$ is an isomorphism. \label{rp4} In particular, if $I$ is in $L$, then $I A^\w$ is the kernel of the canonical projection of $A^\w$ onto $A/I$. 
\ztt
{\bf Proof} \ \ The first claim follows from Proposition VII.4.2 of {\em Homological Algebra} by Cartan and Eilenberg (Princeton University Press, 1956). The second claim follows from Exercise II.2 of the same book in conjunction with Proposition 10.13 of {\em Introduction to Commutative Algebra} by Atiyah and Macdonald (Addison-Wesley, 1969). The third claim is obtained by tensoring the exact sequence 
$$0\to I\to A\to A/I\to0$$
with $A^\w$ over $A$. QED
\section{The Case of an Arbitrary Ring} \label{aq}
(Reminder: by ``ring'' we mean ``commutative ring with 1''.) Let $A$ be a ring, $X$ an indeterminate, and $a$ an element of $A$. For any formal power series 
$$f=\sum_{n\ge0}\ a_n\ (X-a)^n\in A[[X-a]]$$ 
and any nonnegative integer $k$ put 
$$\frac{f^{(k)}}{k!}:=\sum_{n\ge k}\ \bi{n}{k}\ a_n\
(X-a)^{n-k}\ \in A[[X-a]].$$ 
We then get the {\bf Taylor Formula}
$$f=\sum_{n\ge0}\ \frac{f^{(n)}(a)}{n!}\ (X-a)^n.$$

\ii We regard the $A$-algebra $A[[X-a]]$ as an $A[X]$-algebra via the morphism from $A[X]$ to $A[[X-a]]$ which maps $P$ to 
$$\sum_{n\ge0}\ \frac{P^{(n)}(a)}{n!}\ (X-a)^n.$$

\ii For any $f$ in $A[[X-a]]$ and any nonnegative integer $k$ we call $k$-{\bf jet} of $f$ at $a$ the polynomial 
$$\DL_a^k(f)
:=\sum_{n\le k}\ \frac{f^{(n)}(a)}{n!}\ (X-a)^n.$$
Then $\DL_a^k$ induces a ring morphism 
$$\DL_a^k:A[[X-a]]\longrightarrow
\frac{A[X]}{(X-a)^k}\quad.$$

\ii Let $a_1,\dots,a_r$ be in $A$; let $m_1,\dots,m_r$ be positive integers; let $D$ be the product of the $(X-a_i)^{m_i}$; and let $u$ be the natural morphism from $A[X]/(D)$ to the product $B$ of the $A[X]/(X-a_i)^{m_i}$:
$$u:\frac{A[X]}{(D)}\longrightarrow 
B:=\prod_{i=1}^r\ \frac{A[X]}{(X-a_i)^{m_i}}\quad.$$
Assume that $a_i-a_j$ is invertible for all $i\not=j$. 
In particular 
$$D_i(X):=\frac{D(X)}{(X-a_i)^{m_i}}$$
is invertible in $A[[X-a_i]]$. Let $v$ be the morphism 
$$v:\Big(P_i\bmod(X-a_i)^{m_i}\Big)_{i=1}^r\ \mapsto\  
\sum_{i=1}^r\ 
\DL_{a_i}^{m_i-1}\!\!\left(\frac{P_i}{D_i}\right)D_i
\bmod D$$
from $B$ to $A[X]/(D)$. 
\att[Fifth Main Statement] The maps $u$ and $v$ are inverse ring isomorphisms. \label{rp5} Moreover Corollary \ref{e} page \pageref{e} remains true (with the obvious changes of notation). 
\ztt 
{\bf Proof}\ \ The composition $u\circ v$ is obviously the identity of $B$. Since both rings are rank $\deg D$ free $A$-modules, this proves the statement. QED
\section{Universal Remainder}
In this Section we work over some unnamed ring 
(reminder: by ``ring'' we mean ``commutative ring with 1''). The remainder of the euclidean division of $X^r$ by $$(X-a_1)\cdots(X-a_k)$$ is a universal polynomial in $a_1,\dots,a_k,X$ with integral coefficients. To compute this polynomial first recall Newton's interpolation. \bb 

\ii As a general notation set
$$f(u,v):=\frac{f(v)-f(u)}{v-u}\quad,$$
define the {\em Newton interpolation polynomial}  
$$N(f(X);a_1,\dots,a_k;X)$$ 
of an arbitrary polynomial $f(X)$ at $a_1,\dots,a_k$ by

$$N(f(X);a_1,\dots,a_k;X):=f(a_1)+f(a_1,a_2)(X-a_1)$$
$$+f(a_1,a_2,a_3)(X-a_1)(X-a_2)+\cdots$$
$$+f(a_1,\dots,a_k)(X-a_1)\cdots(X-a_{k-1}),$$

and put 
$$g(X):=N(f(X);a_1,\dots,a_k;X),$$
 
$$h(X):=N(f(a_1,X);a_2,\dots,a_k;X).$$ 

\att[Newton Interpolation Theorem] We have $g(a_i)=f(a_i)$ for $i=1,\dots,k$. In particular 
$g(X)$ is the remainder of the euclidean division of $f(X)$ by $(X-a_1)\cdots(X-a_k).$ Moreover $f(a_1,\dots,a_i)$ is symmetric in $a_1,\dots,a_i$.\ztt 
{\bf Proof.} Argue by induction on $k$, the case $k=1$ being easy. Note $$g(X)=f(a_1)+(X-a_1)\ h(X).$$ The equality $g(a_1)=f(a_1)$ is clear. Assume $2\le i\le k$. By induction hypothesis we have $h(a_i)=f(a_1,a_i)$ and thus 
$$g(a_i)=f(a_1)+(a_i-a_1)\ f(a_1,a_i)=f(a_i).$$ 
Since $f(a_1,\dots,a_k)$ is the leading coefficient of $g(X)$, it is symmetric in $a_1,\dots,a_k$. QED \bb

\ii Assume that $$b_i:=\prod_{j\not=i}\ (a_i-a_j)$$ is invertible for all $i$. By Lagrange interpolation we have, for $f(X)=X^r=f_r(X)$,  
$$f_r(a_1,\dots,a_k)=\sum_i\ \frac{a_i^r}{b_i}\quad.$$
The generating function of the sum $s_n=s_n(a_1,\dots,a_k)$ of the degree $n$ monomials in $a_1$, \dots, $a_k$ being 
$$\frac{1}{(1-a_1\,X)\cdots (1-a_k\,X)}
=\sum_i\ \frac{a_i^{k-1}\ b_i^{-1}}{1-a_i\,X}
=\sum_{n,i}\ \frac{a_i^{k-1+n}}{b_i}\ X^n,$$
we get 
$$s_n=\sum_i\ \frac{a_i^{k-1+n}}{b_i}\quad,
\quad f_r(a_1,\dots,a_k)=s_{r-k+1}$$
with the convention $s_n=0$ for $n<0$, and the remainder of the euclidean division of $X^r$ by $(X-a_1)\cdots (X-a_k)$ is 
$$\sum_{i=1}^k\ s_{r-i+1}(a_1,\dots,a_i)\ 
(X-a_1)\cdots(X-a_{i-1})$$
(even when some of the $b_i$ are not invertible). 
\section{An Adjunction} 
Recall that a $\Q${\bf -algebra} is a ring containing $\Q$ (``ring'' meaning ``commutative ring with 1''), and that a {\bf derivation} on a $\Q$-algebra $A$ is a $\Q$-vector space endomorphism $a\mapsto a'$ satisfying $(a b)'=a' b+a b'$ for all $a,b$ in $A$. A {\bf differential} $\Q$-algebra is a $\Q$-algebra equipped with a derivation. We leave it to the reader to define the notions of $\Q$-algebra and of differential $\Q$-algebra {\bf morphism}. To each $\Q$-algebra $A$ is attached the $\Q$-algebra $A[[X]]$ (where $X$ is an indeterminate) equipped with the derivation $\frac{d}{d X}\ $. Let $A$ be a $\Q$-algebra and $B$ a differential $\Q$-algebra. Denote respectively by
$$\mathcal{A}(B,A)\quad\mbox{and}\quad
\mathcal{D}(B,A[[X]])$$
the  $\Q$-vector space of $\Q$-algebra morphisms from $B$ to $A$ and of differential $\Q$-algebra morphisms  from $B$ to $A[[X]]$. The formulas 
$$f(b)=F(b)_0,\quad
F(b)=\sum_{n=0}^\infty\frac{f(b^{(n)})}{n!}\ X^n,$$
where $F(b)_0$ designates the constant term of $F(b)$, set up a bijective $\Q$-linear correspondence between vectors $f$ of $\mathcal{A}(B,A)$ and vectors $F$ of $\mathcal{D}(B,A[[X]])$. Let us check for instance the following point. Let $F$ be in $\mathcal{D}(B,A[[X]])$, let $f$ in $\mathcal{A}(B,A)$ be the constant term of $F$, and let $G$ in $\mathcal{D}(B,A[[X]])$ be the extension of $f$. Let us show $G=F$. We have

$$\begin{array}{lcll}
G(b)&=&\ds\sum_{n=0}^\infty\ \frac{G(b)^{(n)}(0)}{n!}\ X^n&\mbox{by Taylor}\\[2em]
&=&\ds\sum_{n=0}^\infty\ \frac{G(b^{(n)})(0)}{n!}\ X^n&\mbox{because $G(x)'=G(x')$}\\[2em]
&=&\ds\sum_{n=0}^\infty\ \frac{f(b^{(n)})}{n!}\ X^n&\mbox{by definition of $G$}\\[2em]
&=&\ds\sum_{n=0}^\infty\ \frac{F(b^{(n)})(0)}{n!}\ X^n&\mbox{by definition of $f$}\\[2em]
&=&\ds\sum_{n=0}^\infty\ \frac{F(b)^{(n)}(0)}{n!}\ X^n&\mbox{because $F(x')=F(x)'$}\\[2em]
&=&F(b)&\mbox{by Taylor.}
\end{array}$$ 
\section{Wronski} 
Here is a brief summary of Section II.2 in H. Cartan's book {\em Calcul diff\'erentiel}. If $E$ and $F$ are Banach spaces, let $\mathcal{L}(E,F)$ be the Banach space of continuous linear maps from $E$ to $F$. Let $E$ be a Banach space, $I$ a nonempty open interval, and $A$ a continuous map from $I$ to $\mathcal{L}(E,E)$. Consider the ODE's 

\begin{equation}\label{wr1}x'=A(t)\,x,\end{equation} 

where the unknown is a map $x$ from $I$ to $E$, and
\begin{equation}\label{wr2}X'=A(t)\,X,\end{equation}
where the unknown is a map $X$ from $I$ to $F:=\mathcal{L}(E,E)$. Let $t_0$ be in $I$. The solution to (\ref{wr2}) satisfying $X(t_0)=Id_E$ is called the {\bf resolvent} of (\ref{wr1}) and is denoted $t\mapsto R(t,t_0)$.\bb
 
\ii Let $B$ be a continuous map from $I$ to $E$ and $x_0$ a vector of $E$. The solution to 
$$x'=A(t)\,x+B(t),\quad x(t_0)=x_0$$
is 
$$x(t):=R(t,t_0)\,x_0
+\int_{t_0}^t R(t,\tau)\,B(\tau)\ d\tau.$$

\ii Here is an additional comment. For $t$ in $I$ let $\mathcal{A}(t)$ be the left multiplication by $A(t)$. In particular $\mathcal{A}$ is a continuous map from $I$ to $\mathcal{L}(F,F)$. Let $\mathcal{R}$ be its resolvent. Then $\mathcal{R}(t,t_0)$ is the left multiplication by $R(t,t_0)$.\bb

\ii Corollary. If $W$ is a differentiable map from $I$ to $F$ satisfying $W'=A\,W$, then 
$$W(t)=R(t,t_0)\, W(t_0).$$

\ii Corollary to the Corollary. If in addition $W(t_0)$ is invertible, then 
$$R(t,t_0)=W(t)\, W(t_0)^{-1}.$$

\centerline{*}\bb 

[1] For terminological justifications, type {\em didrygaillard} on Google. 
\end{document}